\documentclass[11pt]{article}
\usepackage{geometry}
\usepackage{epsfig}
\usepackage{amsfonts}
\usepackage{amsmath}
\usepackage{amsthm}
\usepackage{graphicx}
\usepackage{fancyhdr}
\usepackage[usenames,dvipsnames]{xcolor}
\usepackage{enumerate}
\usepackage{caption}
\usepackage{subcaption}
\allowdisplaybreaks
 
\newtheorem{theorem}{Theorem}

\newtheorem{proposition}{Proposition}

\newtheorem{lemma}{Lemma}

\newcommand{\RR}{\mathbb{R}}
\newcommand{\VV}{\mathcal{V}}
\newcommand{\FFF}{\mathcal{F}}
\newcommand{\GGG}{\mathcal{G}}
\newcommand{\tX}{\textbf{X}}

\newcommand{\be}{\begin{equation}}
\newcommand{\ee}{\end{equation}}
\newcommand{\bea}{\begin{eqnarray}}
\newcommand{\eea}{\end{eqnarray}}
\newcommand{\bean}{\begin{eqnarray*}}
\newcommand{\eean}{\end{eqnarray*}}
\newcommand{\ben}{\begin{equation*}}
\newcommand{\een}{\end{equation*}}
\newcommand{\ba}{\begin{aligned}}
\newcommand{\ea}{\end{aligned}}

\newcommand{\PP}{\textbf{P}}
\newcommand{\EE}{\textbf{E}}
\newcommand{\ind}{\textbf{1}}

\newcommand{\Nij}{N_n(i,j)}
\newcommand{\NijN}{N_{n+1}(i,j)}
\newcommand{\pij}{p_{ij}}
\newcommand{\nuij}{\nu_n(i,j)}
\newcommand{\Nkl}{N_n(k,l)}

\newcommand{\nukln}{\nu_n(k,l)}

\newcommand{\nuijN}{\nu_{n+1}(i,j)}

\newcommand{\bij}{b^{(i,j)}}
\newcommand{\bijkln}{\bij_{k,l,n}}
\newcommand{\bijklN}{\bij_{k,l,n+1}}
\newcommand{\xiij}{\xi^{(i,j)}}
\newcommand{\xist}{\xi^{(s,t)}}

\newcommand{\DijN}{\Delta_{n+1}(i,j)}
\newcommand{\DklN}{\Delta_{n+1}(k,l)}
\newcommand{\DhfN}{\Delta_{n+1}(h,f)}

\newcommand{\Var}{\mbox{Var}}
\newcommand{\sid}[1]{{\color{black} #1}}
\newcommand{\sidney}[1]{{\color{black} #1}}

\numberwithin{equation}{section}
\numberwithin{theorem}{section}
\numberwithin{corollary}{section}
\numberwithin{proposition}{section}
\numberwithin{lemma}{section}
\numberwithin{figure}{section}
\numberwithin{table}{section}
\usepackage{pdfsync}

\begin{document}
\title{\sc Asymptotic Normality of \sidney{In- and Out-Degree Counts in a} Preferential Attachment Model}
\author{\sc Tiandong Wang and Sidney I. Resnick \thanks{This work was supported by Army MURI grant W911NF-12-1-0385 to Cornell University.}}
\maketitle

\begin{abstract}
  \sidney{Preferential attachment in a directed scale-free graph is widely
  used to model the evolution of social networks.  
Statistical analyses of social networks often relies on node based data
rather than conventional repeated sampling. 
For our directed edge model with preferential attachment, we prove
asymptotic normality of node counts based on a martingale construction
and a martingale central limit theorem. This helps justify estimation
methods based on the statistics of node counts which have specified
in-degree and out-degree.}
\end{abstract}

\section{Introduction}
Preferential attachment for both undirected and directed scale-free
graphs has been introduced in the literature \sid{as a model for the
  growth of} social
networks (cf. \cite{BA99},\cite{BBCR03} and \cite{krapivsky:redner:2001, mori:2002}). 
\sid{Preferential attachment} can \sid{model}
 broader contexts such as the web graph, citation graph, co-author
 graph, etc.  \sid{So more attention} is
now placed upon the directed case where each node \sid{has at least
  two characteristics, namely}
 in- and out-degree. 

Let $\Nij$ be the number of nodes with in-degree $i$ and out-degree
$j$ in a simplified directed preferential attachment model at $n$th
step of the growth of the network. \cite{BBCR03}
showed that $\Nij/n\rightarrow p_{ij}$ for fixed $i$ and $j$, and
\sid{provided an} explicit form of $(p_{ij})$.
Furthermore, we know that the limiting degree sequence $(p_{ij})$ has
both marginally and jointly regularly varying tails (cf. \cite{BBCR03}, \cite{RS15} and 
\cite{SRTDWW16}). However, \sid{what remains an
  open issue for these models is} rigorous justification of
\sid{methods of}  statistical analyses \sid{based on data from}
 social networks. \sid{Therefore, the object of this paper}
 is to examine the asymptotic normality of $\Nij$ \sid{with the idea that
   this asymptotic normality  can} justify  statistical estimation methods in
 practice. Using the martingale central limit theorem, we will show \sid{this}
 asymptotic normality of $\sqrt{n}(\Nij/n-p_{ij})$ for fixed
 $(i,j)$ \sid{as well as jointly over $(i,j)$}. Hence, we 
conclude that the empirical estimator \sid{$\Nij/n$} is
 consistent and asymptotically normal. \sid{We will explore more formal 
 statistical inference that relies on node based
 data and the asymptotic normality elsewhere and give examples of analyses.}

The  directed preferential attachment model that we study 
is outlined in Section~\ref{s2} and our main results on normality are summarised in
Section~\ref{s3}. Proofs are collected in Section~\ref{s4}. 

\section{Model}\label{s2}
\sid{We somewhat simplify the model used in
\cite{BBCR03, RS15,SRTDWW16}.}
At each step of the construction, a node is added;
 we exclude the possibility of adding only a new edge between
existing nodes. The model evolves according to the following
dynamics. 
Choose strictly positive parameters $\alpha,\gamma,\lambda$ and $\mu$ such that $\alpha+\gamma=1$, and we assume in addition that $\alpha,\gamma<1$ to avoid trivial cases.

We initiate the algorithm with a simple case: a graph $G_1$ with one single node (labeled 1) with a self-loop so that both its in and out degrees are 1, denoted by $D_1(1)=(1,1)$.
At stage $n$, we have a directed random graph $G_n=(V_n, E_n)$. If a
node $v$ is from $V_n$, use $D_{in}(v)$ and $D_{out}(v)$ to denote its
in and out degree respectively \sid{(dependence on $n$ is suppressed)} and write
$D_n(v)=(D_{in}(v),D_{out}(v))$. Then $G_{n+1}$ is obtained from $G_n$
as follows. 
\begin{enumerate}
\item[(i)] With probability $\alpha$ a new node $w$ is born and we add an edge leading from $w$ to an existing node $v\in V_n$. The existing node $v$ is chosen with probability according to its in-degree:
\be\label{in}\PP(v\in V_n\mbox{ is chosen})=\frac{D_{in}(v)+\lambda}{(1+\lambda)n}.\ee
\item[(ii)] With probability $\gamma$ a new node $w$ is born and we add an edge leading from an existing node $v\in V_n$ to $w$. The existing node $v$ is chosen with probability according to its out-degree:
\be\label{out}\PP(v\in V_n\mbox{ is
  chosen})=\frac{D_{out}(v)+\mu}{(1+\mu)n}.\ee 
\end{enumerate}
\sid{The construction makes} $G_n$  a directed graph with $n$ nodes
(i.e. $V_n=\{1,2,\ldots,n\}$) and $n-1$ edges; \sid{the self-loop in $G_1$
is not counted as an edge}. Note that  
\[\sum_{v\in V_n}D_{in}(v)= \sum_{v\in V_n}D_{out}(v)=n,\]
so the attachment probabilities in \eqref{in} and \eqref{out} add to 1.

\section{Results}\label{s3}
For $i,j\ge 0$, let $\Nij$ denote the number of nodes with in-degree $i$ and out-degree $j$ in $G_n$, i.e.
\[\Nij=\sum_{v\in V_n}\ind_{\{D_n(v)=(i,j)\}}\quad \sid{(n\geq 1)},\]
\sid{and}  set $\nuij=\EE(\Nij)$. The following lemma elaborates part of
the results of Theorem~3.2 in \cite{BBCR03}, which
implies that for each $i$ and $j$ there are non-random constants
$(p_{ij})$ such that 
\be\label{convij}
\frac{\Nij}{n}\rightarrow \pij \mbox{ a.s. as } n\rightarrow\infty,
\ee
Clearly, $p_{00}=0$. We also take $\Nij$ and $p_{ij}$ to be zero if
either $i$ or $j$ is $-1$. \sid{The explicit form of the} limiting degree
distribution  $(p_{ij})$ is given in \cite{BBCR03}. 
\begin{lemma}\label{degseq}
For each $i,j=0,1,2,\ldots$, we have for $C>4$,
\be\label{step1}
|\nuij-n\pij|\le C, \mbox{ for }\forall n\ge 1, \quad\forall (i,j),
\ee and
\be\label{concentration}
\PP\left(\bigvee_{(i,j)}\left|\frac{\Nij}{n}-p_{ij}\right|\ge C\sqrt{\frac{\log n}{n}}\right)=o(1), \mbox{ as } n\rightarrow\infty,
\ee
where \sid{the} $\pij$ satisfy
\be\label{pij}
\pij=\alpha \ind_{\{(i,j)=(0,1)\}}+\gamma \ind_{\{(i,j)=(1,0)\}}+c_1(i-1+\lambda)p_{i-1,j}+c_2(j-1+\mu)p_{i,j-1}-\delta_{ij}\pij.
\ee
Here we have
\[c_1=\frac{\alpha}{1+\lambda},\quad c_2=\frac{\gamma}{1+\mu},\quad \delta_{ij}=c_1(i+\lambda)+c_2(j+\mu).\]
\end{lemma}
\sid{As a stochastic process in $(i,j)$, the proportion
of nodes with in-degree $i$ and out-degree $j$ converges in
distribution after centering and scaling to a centered Gaussian process.}
Asymptotic
  normality relies on a standard multivariate martingale central limit theorem
(cf. Proposition~2.2 outlined in~\cite{counts}; \sid{a statement is
  given in Proposition \ref{mgclt} in  Section \ref{mishmash}} and see
also
\cite{kuchler:sorensen:1999,hall:heyde:1980,
hubalek:posedel:2013, crimaldi:pratelli:2005} and \cite[Chapter
8]{durrett:2010a}).
For our
problem, the normality results are 
summarized in the next theorem. 

\begin{theorem}\label{main}  
\sid{Fix positive integers $I,O$. In the normality statement, matrices
  $K_{IO}$ and $\Sigma_{IO}$} are 
are specified in \eqref{KIO} and \eqref{matrixAsyVar} respectively.  
 Provided that $K_{IO}$ is invertible, we have
\be\label{thm}\left(\sqrt{n}\left(\frac{N_n(i,j)}{n}-p_{ij}\right):\,  0\le i\le I,\, 0\le j\le O\right)\Rightarrow N(0,K^{-1}_{IO}\Sigma_{IO}K^{-T}_{IO}).\ee
\end{theorem}
 
\section{Proofs}\label{s4}

%\begin{lemma}\label{degseq}
%For each $i,j=0,1,2,\ldots$, we have
%\be\label{convij}
%\frac{\Nij}{n}\rightarrow \pij \mbox{ a.s. as } n\rightarrow\infty,
%\ee
%where $\pij$ satisfy
%\be\label{pij}
%\pij=\alpha \ind_{\{(i,j)=(0,1)\}}+\gamma \ind_{\{(i,j)=(1,0)\}}+c_1(i-1+\lambda)p_{i-1,j}+c_2(j-1+\mu)p_{i,j-1}-\delta_{ij}\pij.
%\ee
%Here we have
%\[c_1=\frac{\alpha}{1+\lambda},\quad c_2=\frac{\gamma}{1+\mu},\quad \delta_{ij}=c_1(i+\lambda)+c_2(j+%\mu).\]
%\end{lemma}
%\begin{proof}
\subsection{Proof of Lemma~\ref{degseq}} 
By the construction of our model, at the initial stage we have
$N_1(1,1)=1$, $N_1(i,j)=0$ for $(i,j)\neq (1,1)$. \sid{Let $\FFF_n$ be
  the $\sigma$-field of information accumulated by watching the graph
  grow until stage $n$. We have}
%\begin{eqnarray*}
\begin{align}
\EE(N_{n+1}(i,j)|\FFF_n)
=& \Nij+\EE(N_{n+1}(i,j)-\Nij|\FFF_n)\nonumber \\
=& \Nij+\alpha\ind_{\{(i,j)=(0,1)\}}+\gamma \ind_{\{(i,j)=(1,0)\}}\nonumber \\
&\:+\PP(\mbox{a new edge from $n+1$ to $v\in V_n$;}
  D_n(v)=(i-1,j)|\FFF_n)\nonumber  \\
&\:+\PP(\mbox{a new edge from $v\in V_n$ to $n+1$;}
  D_n(v)=(i,j-1)|\FFF_n)\nonumber \\ 
&\:-\PP(\mbox{a new edge from $n+1$ to $v\in V_n$;}
  D_n(v)=(i,j)|\FFF_n)\nonumber  \\
&\:-\PP(\mbox{a new edge from $v\in V_n$ to $n+1$;} D_n(v)=(i,j)|\FFF_n)\nonumber \\
=& \Nij+\alpha\ind_{\{(i,j)=(0,1)\}}+\gamma \ind_{\{(i,j)=(1,0)\}}\nonumber \\
&\:+
  c_1(i-1+\lambda)\frac{N_n(i-1,j)}{n}+c_2(j-1+\mu)\frac{N_n(i,j-1)}{n}\nonumber 
  \\ 
&\:- (c_1(i+\lambda))+c_2(j+\mu))\frac{\Nij}{n}. \label{siddef}
%\end{eqnarray*}
\end{align}
\sid{Taking expectations and recalling that}  $\nuij:=\EE(\Nij)$, we get
\bea
\nu_{n+1}(i,j)&=&\alpha\ind_{\{(i,j)=(0,1)\}}+\gamma \ind_{\{(i,j)=(1,0)\}}\nonumber\\
&&\:+(1-\frac{\delta_{ij}}{n})\nuij+
\frac{c_1(i-1+\lambda)}{n}\nu_n(i-1,j)+\frac{c_2(j-1+\mu)}{n}\nu_n(i,j-1).\label{nudef}
\eea

We first show that \eqref{step1} holds for some constant $C\ge 1$.
Let $\varepsilon_n(i,j)=\nuij-n\pij$, then
\be\label{base}|\varepsilon_1(1,1)|=|1-p_{11}|\le C,\quad |\varepsilon_1(i,j)|=|0-p_{ij}|\le C\mbox{ for }(i,j)\neq (1,1).\ee
Also, for $n\ge 1$ 
\[\varepsilon_{n+1}(0,1)=(1-\frac{\delta_{01}}{n})\varepsilon_n(0,1),\quad
 \varepsilon_{n+1}(1,0)=(1-\frac{\delta_{10}}{n})\varepsilon_n(1,0),\]
  and further for $(i,j)\notin\{(0,1),(1,0)\}$:
\bean
\varepsilon_{n+1}(i,j)&=& \nu_{n+1}(i,j)-(n+1)\pij\\
&=& (1-\frac{\delta_{ij}}{n})\varepsilon_n(i,j)+\frac{c_1(i-1+\lambda)}{n}\varepsilon(i-1,j)+\frac{c_2(j-1+\mu)}{n}\varepsilon(i,j-1).
\eean
Then \eqref{step1} is true for $(i,j)=(0,1)$ by simpl\sid{e} induction on $n$: if $|\varepsilon_n(0,1)|\le C$, then 
\[|\varepsilon_{n+1}(0,1)|\le \left|1-\frac{\delta_{ij}}{n}\right|C\le C.\] Similar arguments give that \eqref{step1} also holds for $(i,j)=(1,0)$. For $(i,j)\notin\{(0,1),(1,0)\}$,
our induction assumption $\bigvee_{(i,j)}|\varepsilon_n(i,j)|\le C$ (which is true for $n=1$ by \eqref{base}) gives that
\bean
|\varepsilon_{n+1}(i,j)|&=& \left|(1-\frac{\delta_{ij}}{n})\varepsilon_n(i,j)+\frac{c_1(i-1+\lambda)}{n}\varepsilon_n(i-1,j)+\frac{c_2(j-1+\mu)}{n}\varepsilon_n(i,j-1)\right|\\
&\le & \left|1-\frac{\delta_{ij}}{n}+\frac{c_1(i-1+\lambda)}{n}+\frac{c_2(j-1+\mu)}{n}\right|C\\
&=& \left|1-\frac{c_1+c_2}{n}\right|C\le C,
\eean
by noting that $c_1+c_2\le \alpha+\gamma=1$.
Hence,
\[\bigvee_{(i,j)}|\varepsilon_{n+1}(i,j)|=|\varepsilon_{n+1}(0,1)|\vee|\varepsilon_{n+1}(1,0)|\vee\bigvee_{(i,j)\notin\{(0,1),(1,0)\}}|\varepsilon_{n+1}(i,j)|\le C.\]
\sid{This verifies that  \eqref{step1} holds.}

Next fix $(i,j)$ \sid{and $n$} and define the uniformly integrable
martingale 
$$Y_m(i,j)=\EE(\Nij|\FFF_m),\quad m=0,1,\ldots,n$$ \sid{with
difference sequence }
\[|d_m(i,j)|=|Y_m(i,j)-Y_{m-1}(i,j)|.\]
By~\cite{BBCR03}, given $\FFF_m$, determining $G_n$ only requires the identification 
of which old vertices are involved at each stage, and there are at most $2n$ such choices. 
Under proper redistribution, changing one of these choices (say from node $u$ to node $v$)
 will only alter the degrees of $u$ and $v$ in the final graph. Hence,
 \[|d_m(i,j)|\le 2, \quad m=0,1,\ldots,n.\]

Also, $Y_0(i,j)=\EE(\Nij|\FFF_0)=\nuij$, and 
\[\Nij-\nuij=\sum_{k=1}^n d_k(i,j).\]
Then by the Azuma-Hoeffding's inequality \cite{azuma:1967}
\[\PP(|\Nij-\nuij|\ge C\sqrt{n\log n})\le 2\exp\left(-\frac{C^2 n\log n}{2n\cdot 4}\right)=\frac{2}{n^{C^2/8}}.\]
Therefore,
\bean
\lefteqn{\PP\left(\bigvee_{(i,j)}|\Nij-\nuij|\ge C\sqrt{n\log n}\right)}\\
&\le & n^2\bigvee_{(i,j)}\PP(|\Nij-\nuij|\ge C\sqrt{n\log n})\\
&\le & 2n^{-(C^2/8-2)}
\eean
In other words, for $C>4$,
\be\label{step2}
\PP\left(\bigvee_{(i,j)}|\Nij-\nuij|\ge C\sqrt{n\log n}\right)=o(1).
\ee
Now \eqref{step1} and \eqref{step2} together imply that  
%\marginpar{\sid{Why -1?}}
\[\PP\left(\bigvee_{(i,j)}\left|\frac{\Nij}{n}-p_{ij}\right|\ge
  \frac{C}{n}(\sqrt{n\log n}-1)\right)=o(1),\] 
\sid{and this gives} \eqref{concentration}
for $C>4$.

\subsection{Proof of Theorem~\ref{main}}\label{mishmash}
\subsubsection{Preliminary: Martingale central limit theorem.}\label{prelim}
We  use a multivariate martingale central limit theorem to prove 
Theorem~\ref{main}. \sid{We first state the version that we need. See
Proposition~2.2 in \cite{RS15} and also
\cite{kuchler:sorensen:1999,hall:heyde:1980, 
hubalek:posedel:2013, crimaldi:pratelli:2005} and \cite[Chapter
8]{durrett:2010a}).}

\begin{proposition}\label{mgclt}
Let $\{\tX_{n,m},\GGG_{n,m},1\le m\le n\}$,
$\tX_{n,m}=(X_{n,m,1},\ldots,X_{n,m,d})^T$, be a d-dimensional
square-integrable \sid{martingale} difference array. Consider the
$d\times d$ nonnegative definite random matrices 
\[G_{n,m}=\left(\EE(X_{n,m,i}X_{n,m,j}|\GGG_{n,m-1}),i,j=1,2,\ldots, d\right),\quad V_n=\sum_{m=1}^n G_{n,m},\]
and suppose $(A_n)$ is a sequence of $l\times d$ matrices with a bounded supremum norm. Assume that
 \begin{enumerate}
 \item[(i)] $A_nV_n A_n^T\stackrel{P}{\rightarrow}\Sigma$ as $n\rightarrow\infty$ for some deterministic (automatically nonnegatively definite) matrix $\Sigma$.
 \item[(ii)] $\sum_{m\le n}\EE(X^2_{n,m,i}1_{\{|X_{n,m,i}|>\epsilon\}}|\GGG_{n,m-1})\stackrel{P}{\rightarrow}0$ as $n\rightarrow\infty$ for all $i=1,2,\ldots,d$ and $\epsilon>0$.
 \end{enumerate}
 Then in $\RR^l$, as $n\rightarrow\infty$
 \be\label{clt}
 \sum_{m=1}^n A_n\tX_{n,m}\Rightarrow \tX,
 \ee 
a centered l-dimensional Gaussian vector with covariance matrix $\Sigma$.
\end{proposition}

\subsubsection{The martingale.}\label{themg}
We start with constructing a martingale for fixed $i$ and $j$. Suppose that our martingale takes the form
\be\label{mg}
M_n(i,j)=\sum_{l=0}^j\sum_{k=0}^i b_{k,l,n}^{(i,j)}(N_n(k,l)-\nu_n(k,l)),
\ee
where $b_{k,l,n}^{(i,j)}$ are some non-random constants. \sid{We
investigate what properties $b_{k,l,n}^{(i,j)}$ must satisfy in order
that $M_n(i,j)$ is a martingale in $n$.}

By the model assumptions in Section~\ref{s2}, $n\ge i\vee j$ since 
at each stage we
can only increase either the in or out-degree of one particular node
by 1. Therefore, with probability 1, $\Nij=0$ for
$n<i\vee j$. Also for $n=1$, almost surely, $N_1(1,1)=1$ and
$N_1(i,j)=0$ for all other $(i,j)\neq(1,1)$. Hence,  
\[N_1(i,j)-\nu_1(i,j)=0 \quad\mbox{a.s. for all }(i,j),\]
and values of $b_{k,l,1}^{(i,j)}$ will not affect the form of $M_1(i,j)$. For simplicity of calculations, we set $b_{k,l,1}^{(i,j)}=1$ for all $(k,l)$.

Using \eqref{siddef} and \eqref{nudef}, we see that in order to make $M_n(i,j)$ a martingale, we must have
\bean
\EE(M_{n+1}(i,j)|\FFF_n)&=&\sum_{l=0}^j\sum_{k=0}^i\bijklN\left[(1-\frac{\delta_{kl}}{n})(\Nkl-\nukln)\right.\\
&&\:\left.+\frac{c_1(k-1+\lambda)}{n}(N_n(k-1,l)-\nu_n(k-1,l))\right.\\
&&\:\left.+\frac{c_2(l-1+\mu)}{n}(N_n(k,l-1)-\nu_n(k,l-1))\right]\\
&=& \sum_{l=0}^j\sum_{k=0}^i b_{k,l,n}^{(i,j)}(N_n(k,l)-\nu_n(k,l))=M_n(i,j),
\eean
where the last equality follows from the martingale assumption. Thus, $\bijkln, 0\le k\le i, 0\le l\le j$, must satisfy the following recursions:
\begin{align}
\bij_{i,j,n+1}\left(1-\frac{\delta_{ij}}{n}\right)&=\bij_{i,j,n}\label{rec1}\\
\bij_{k,j,n+1}\left(1-\frac{\delta_{kj}}{n}\right)+\bij_{k+1,j,n+1}\frac{c_1(k+\lambda)}{n}&=\bij_{k,j,n},\quad 0\le k\le i-1 \label{rec2}\\
\bij_{i,l,n+1}\left(1-\frac{\delta_{il}}{n}\right)+\bij_{i,l+1,n+1}\frac{c_2(l+\mu)}{n}&=\bij_{i,l,n},\quad 0\le l\le j-1 \label{rec3}\\
\bijklN\left(1-\frac{\delta_{kl}}{n}\right)+\bij_{k+1,l,n+1}\frac{c_1(k+\lambda)}{n}+\bij_{k,l+1,n+1}\frac{c_2(l+\mu)}{n}
&=\bij_{k,l,n},\nonumber\\
&\quad 0\le k\le i-1, 0\le l\le j-1.\label{rec4}
\end{align}
Solving \eqref{rec1} gives
\begin{equation}\label{prod}
\bij_{i,j,n+1}=\prod_{m=1}^n\left(1-\frac{\delta_{ij}}{m}\right)^{-1}.
\end{equation}
Also, \eqref{rec2} yields
\begin{align}
\bij_{k,j,n+1}=&\left(1-\frac{\delta_{kj}}{n}\right)^{-1}\bij_{k,j,n}-\left(1-\frac{\delta_{kj}}{n}\right)^{-1}\bij_{k+1,j,n+1}\frac{c_1(k+\lambda)}{n}
  \nonumber\\
=&
   \left(1-\frac{\delta_{kj}}{n}\right)^{-1}\left[\left(1-\frac{\delta_{kj}}{n-1}\right)^{-1}\bij_{k,j,n-1}-\left(1-\frac{\delta_{kj}}{n-1}\right)^{-1}\bij_{k+1,j,n}\frac{c_1(k+\lambda)}{n-1}\right]
   \nonumber \\   
&\quad
  \:-\left(1-\frac{\delta_{kj}}{n}\right)^{-1}\bij_{k+1,j,n+1}\frac{c_1(k+\lambda)}{n}
 \nonumber\\
=&\ldots =
   \prod_{m=1}^n\left(1-\frac{\delta_{kj}}{m}\right)^{-1}-\sum_{m=1}^n
   \bij_{k+1,j,m+1}\frac{c_1(k+\lambda)}{m}\prod_{d=m}^n\left(1-\frac{\delta_{kj}}{d}\right)^{-1}. \label{prod1}
\end{align}
Similarly, we can obtain from \eqref{rec3} \sid{and}  \eqref{rec4} that
\[\bij_{i,l,n+1}=\prod_{m=1}^n\left(1-\frac{\delta_{il}}{m}\right)^{-1}-\sum_{m=1}^n \bij_{i,l+1,m+1}\frac{c_2(l+\mu)}{m}\prod_{d=m}^n\left(1-\frac{\delta_{il}}{d}\right)^{-1},\]
and that
\begin{align*}
\bijklN=\prod_{m=1}^n\left(1-\frac{\delta_{kl}}{m}\right)^{-1}&-\sum_{m=1}^n \bij_{k+1,l,m+1}\frac{c_1(k+\lambda)}{m}\prod_{d=m}^n\left(1-\frac{\delta_{kl}}{d}\right)^{-1}\\
&-\sum_{m=1}^n \bij_{k,l+1,m+1}\frac{c_2(l+\mu)}{m}\prod_{d=m}^n\left(1-\frac{\delta_{kl}}{d}\right)^{-1}.
\end{align*}

\subsubsection{Properties of the coefficients ${\bijklN}$}
\sid{For the calculation of the asymptotic form of conditional
  covariances of martingale differences, we will need} the asymptotic forms of the ratio
\sid{${\bijklN}/{\bij_{i,j,n+1}}$} for all $k\le i, l\le j$, as
$n\rightarrow\infty$ and we set  
\sid{\begin{equation}\label{xi}
\xi_{kl}^{(i,j)}:=\lim_{n\rightarrow\infty}\frac{\bij_{k,l,n+1}}{\bij_{i,j,n+1}},\quad 
k=0,1,\ldots ,i \text{ and } l=0,1,\ldots ,j.
\end{equation}}
 \sid{We begin with the case  $l=j$}. \sid{Using \eqref{prod} and \eqref{prod1}} we know that for $0\le k\le i-1$, 
\be\label{bratio}
\frac{\bij_{k,j,n+1}}{\bij_{i,j,n+1}}=\prod_{m=1}^n\frac{\left(1-\frac{\delta_{kj}}{m}\right)^{-1}}{\left(1-\frac{\delta_{ij}}{m}\right)^{-1}}-\sum_{m=1}^n\frac{\bij_{k+1,j,m+1}}{m}c_1(k+\lambda)\frac{\prod_{d=m}^n\left(1-\frac{\delta_{kj}}{d}\right)^{-1}}{\prod_{m=1}^n\left(1-\frac{\delta_{ij}}{m}\right)^{-1}},\ee 
\sid{and  from \eqref{bratio},} we claim that
\be\label{limbkj}
\frac{\bij_{k,j,n+1}}{\bij_{i,j,n+1}}\rightarrow \sid{\xi_{kj}^{(i,j)}}= (-1)^{i-k}\prod_{d=k}^{i-1}\left(\frac{\lambda+d}{i-d}\right).
\ee
\sid{For the first term on the right of \eqref{bratio} we have} by Stirling's formula,
\begin{equation}\label{littleoh}
\prod_{m=1}^n\frac{\left(1-\frac{\delta_{kj}}{m}\right)^{-1}}{\left(1-\frac{\delta_{ij}}{m}\right)^{-1}}=\prod_{m=1}^n\frac{m-\delta_{ij}}{m-\delta_{kj}}\sim\frac{\Gamma(1-\delta_{ij})}{\Gamma(1-\delta_{kj})}n^{-(i-k)c_1}\rightarrow
0,\mbox{ as }n\rightarrow\infty,
\end{equation}
because $i-k\ge 1$.
\sid{Hence proving \eqref{limbkj} requires showing}
\be\label{claimbij}
-\sum_{m=1}^n\frac{\bij_{k+1,j,m+1}}{\bij_{i,j,n+1}}\frac{c_1(k+\lambda)}{m}\prod_{d=m}^n\left(1-\frac{\delta_{kj}}{d}\right)^{-1}\rightarrow (-1)^{i-k}\prod_{d=k}^{i-1}\left(\frac{\lambda+d}{i-d}\right)
\ee
\sid{and w}e prove this by induction on $k<i$. For $k=i-1$, \sid{using
\eqref{prod},} we have
\begin{align}
-\sum_{m=1}^n & \frac{\bij_{i,j,m+1}}{m}c_1(i-1+\lambda)\frac{\prod_{d=m}^n\left(1-\frac{\delta_{i-1,j}}{d}\right)^{-1}}{\prod_{d=1}^n\left(1-\frac{\delta_{ij}}{d}\right)^{-1}}\nonumber\\
=& -c_1(i-1+\lambda)\sum_{m=1}^n\frac{1}{m}\frac{\prod_{d=m}^n\left(1-\frac{\delta_{i-1,j}}{d}\right)^{-1}}{\prod_{d=m+1}^n\left(1-\frac{\delta_{ij}}{d}\right)^{-1}}\nonumber\\
=& -c_1(i-1+\lambda)\sum_{m=1}^n\frac{\Gamma(n+1-\delta_{ij})/\Gamma(m+1-\delta_{ij})}{\Gamma(n+1-\delta_{i-1,j})/\Gamma(m-\delta_{i-1,j})}\nonumber\\
=&-c_1(i-1+\lambda)\frac{\Gamma(n+1-\delta_{ij})}{\Gamma(n+1-\delta_{i-1,
   j})}\sum_{m=1}^n g(m),  \label{xibase}
\end{align}
\sid{where }
$$
g(m)= 
\frac{\Gamma (m-\delta_{i-1,j})}{\Gamma (m+1-\delta_{ij})}.
$$
Stirling's formula gives as $n\to \infty$
$$
\frac{\Gamma(n+1-\delta_{ij})}{\Gamma(n+1-\delta_{i-1,
   j})} \sim n^{\delta_{i-1,j}-\delta_{ij}} =n^{-c_1}$$ 
and also
$$g(n) \sim n^{c_1-1},\quad (n\to\infty).$$
So the function $g(n)$ is regularly varying and hence by Karamata's
theorem on integration (see, for example, \cite{resnickbook:2007}), we
have
$$\sum_{m=1}^n g(m) \sim n^{c_1}/c_1$$ 
and thus \eqref{xibase} is asymptotic to
$$ -c_1(i-1+\lambda)n^{-c_1} n^{c_1}/c_1=-(i-1+\lambda).$$
This verifies the base case for \eqref{claimbij} and thus \eqref{limbkj} is also true when $k=i-1$. 

\sid{For the next step in the induction argument, we s}uppose that
\eqref{claimbij} holds for $k+1$. Then because of \eqref{littleoh},
\eqref{limbkj} holds for $k$. \sid{We then evaluate the left side of
\eqref{claimbij} with $k+1$ replaced by $k$.}  \sid{Using
\eqref{prod},
$\Gamma(t+1)=t\Gamma(t)$  and 
calculations similar to what was just done, we get}
\begin{align}
-&\sum_{m=1}^n \frac{\bij_{k,j,m+1}}{\bij_{i,j,n+1}}\frac{c_1(k-1+\lambda)}{m}\prod_{d=m}^n\left(1-\frac{\delta_{k-1,j}}{d}\right)^{-1}\nonumber\\
&
= -\sum_{m=1}^n \frac{\bij_{k,j,m+1}}{\bij_{i,j,m+1}}\frac{c_1(k-1+\lambda)}{m}\frac{\prod_{d=m}^n\left(1-\frac{\delta_{k-1,j}}{d}\right)^{-1}}{\prod_{d=m+1}^n\left(1-\frac{\delta_{ij}}{d}\right)^{-1}}\nonumber\\
&= -c_1(k-1+\lambda)\frac{\Gamma(n+1-\delta_{ij})}{\Gamma(n+1-\delta_{k-1, j})}\sum_{m=1}^n\left\{
\frac{\bij_{k,j,m+1}}{\bij_{i,j,m+1}}\frac{1}{m-\delta_{k-1,j}}
\frac{\Gamma(1-\delta_{k-1,j})}{\Gamma(1-\delta_{ij})}
\prod_{d=1}^m\frac{d-\delta_{k-1,j}}{d-\delta_{ij}}\right\}\nonumber\\
&= \sid{-c_1(k-1+\lambda)\frac{\Gamma(n+1-\delta_{ij})}{\Gamma(n+1-\delta_{k-1, j})}\sum_{m=1}^nh(m),}\label{xiinduc}
\end{align}
%\lefteqn{-\sum_{m=1}^n \frac{\bij_{k,j,m+1}}{\bij_{i,j,n+1}}\frac{c_1(k-1+\lambda)}{m}\prod_{d=m}^n\left(1-\frac{\delta_{k-1,j}}{d}\right)^{-1}}\\
%&\sim & (-1)^{i-k}\prod_{d=k}^{i-1}\left(\frac{\lambda+d}{i-d}\right)c_1(k-1+\lambda)\left(-\sum_{m=1}^n\frac{1}{m}\prod_{d=m}^n\left(\frac{d-\delta_{ij}}{d-\delta_{k-1,j}}\right)\right)\\
%& \sim & (-1)^{i-k+1}\prod_{d=k}^{i-1}\left(\frac{\lambda+d}{i-d}\right)(k-1+\lambda)\int_0^1 c_1 x^{(i-k+1)c_1-1}\dd x\\
%&=& (-1)^{i-k+1}\prod_{d=k-1}^{i-1}\left(\frac{\lambda+d}{i-d}\right).
\sid{where}
\begin{align*}
h(m)&= \frac{\bij_{k,j,m+1}}{\bij_{i,j,m+1}}\frac{1}{m-\delta_{k-1,j}}
\frac{\Gamma(1-\delta_{k-1,j})}{\Gamma(1-\delta_{ij})}
\prod_{d=1}^m\frac{d-\delta_{k-1,j}}{d-\delta_{ij}}\\
&=
  \frac{\bij_{k,j,m+1}}{\bij_{i,j,m+1}}\times\frac{\Gamma(m+1-\delta_{k-1,j})}
{\sid{\Gamma ( m+1-\delta_{ij})}}\frac{1}{m-\delta_{k-1,j}}. \end{align*}
\sid{Since the induction assumption means that
  \eqref{limbkj} holds for $k$, we have, as $m\to\infty,$ that $h$ is
  regularly varying with index $(i-k+1)c_1-1$,}
$$h(m) {\sim} m^{(i-k+1)c_1-1}(-1)^{i-k}\prod_{d=k}^{i-1}\left(\frac{\lambda+d}{i-d}\right).
$$
Again, using Karamata's theorem, we have from \eqref{xiinduc}:
\begin{align*}
-c_1(k-1+\lambda)\frac{\Gamma(n+1-\delta_{ij})}{\Gamma(n+1-\delta_{k-1,
  j})}\sum_{m=1}^nh(m) 
&\rightarrow
  (-1)^{i-k+1}\frac{k-1+\lambda}{i-k+1}\prod_{d=k}^{i-1}\left(\frac{\lambda+d}{i-d}\right) \\
&= (-1)^{i-k+1}\prod_{d=k-1}^{i-1}\left(\frac{\lambda+d}{i-d}\right).
\end{align*}
Hence \eqref{claimbij} holds for all $k=0,1,\ldots, i-1$. 
With $\xi_{kj}^{(i,j)}$ defined in \eqref{xi}, we have verified \eqref{limbkj}.

Similarly, as $n\rightarrow\infty$,
\begin{align}
\frac{\bij_{i,l,n+1}}{\bij_{i,j,n+1}} & \rightarrow (-1)^{j-l}\prod_{r=k}^{j-1}\left(\frac{\mu+r}{j-r}\right) =:\xi_{il}^{(i,j)}\label{limbil},\\
\frac{\bij_{k,l,n+1}}{\bij_{i,j,n+1}} &\sim -\left[\frac{c_1(k+\lambda)}{\delta_{ij}-\delta_{kl}}\times\frac{\bij_{k+1,l,n+1}}{\bij_{i,j,n+1}}+\frac{c_2(l+\mu)}{\delta_{ij}-\delta_{kl}}\times\frac{\bij_{k,l+1,n+1}}{\bij_{i,j,n+1}}\right]\nonumber\\
&\rightarrow
  -\left[\frac{c_1(k+\lambda)}{\delta_{ij}-\delta_{kl}}\times\xi_{k+1,l}^{(i,j)}+\frac{c_2(l+\mu)}{\delta_{ij}-\delta_{kl}}\times\xi_{k,l+1}^{(i,j)}\right]=:\xi_{kl}^{(i,j)},\quad
  0\le k\le i-1,\;0\le l\le j-1.\label{limbkl} 
\end{align}
%Define $\xi_{kl}^{(i,j)}:=\lim_{n\rightarrow\infty}\frac{\bij_{k,l,n+1}}{\bij_{i,j,n+1}}$ for all $k=0,1,\ldots i$ and $l=0,1,\ldots j$, then according to \eqref{limbkj}, \eqref{limbil} and \eqref{limbkl}, it satisfies the following:
%\begin{align}
%\xi_{ij}^{(i,j)}&=1, \xi_{00}^{(i,j)}=0,\label{xispe}\\
%\xi_{kj}^{(i,j)}&=(-1)^{i-k}\prod_{d=k}^{i-1}\left(\frac{\lambda+d}{i-d}\right),\quad 0\le k\le i-1,\label{xikj}\\
%\xi_{il}^{(i,j)}&=(-1)^{j-l}\prod_{r=k}^{j-1}\left(\frac{\mu+r}{j-r}\right),\quad 0\le l\le j-1, \label{xibil}\\
%\xi_{kl}^{(i,j)}&=-\left[\frac{c_1(k+\lambda)}{\delta_{ij}-\delta_{kl}}\times\xi_{k+1,l}^{(i,j)}+\frac{c_2(l+\mu)}{\delta_{ij}-\delta_{kl}}\times%\xi_{k,l+1}^{(i,j)}\right],\quad 0\le k\le i-1\mbox{ and }0\le l\le j-1.\label{xibkl}
%\end{align}
We set $\xi_{ij}^{(i,j)}=1, \xi_{00}^{(i,j)}=0$, and note that $\xi_{kl}^{(i,j)}=0$ if either $k>i$ or $l>j$.

\subsubsection{Martingale differences.}\label{secmgdiff} 
Now we are ready to consider the martingale difference:
\begin{align}\label{mgdiff}
M&_{n+1}  (i,j)-M_{n}(i,j)\nonumber\\
&=\sum_{l=0}^j\sum_{k=0}^i (\bij_{k,l,n+1}N_{n+1}(k,l)-\bijkln N_{n}(k,l))
-\sum_{l=0}^j\sum_{k=0}^i (\bij_{k,l,n+1}\nu_{n+1}(k,l)-\bijkln\nu_{n}(k,l)).
\end{align}
\sid{Consider the second double sum on the right side} of
\eqref{mgdiff}. Recall that $\nuij$ satisfies the recursion in
\eqref{nudef}, and this together with the properties of $\bijkln$
in \eqref{rec1}--\eqref{rec4}) give
\begin{align*}
\sum_{l=0}^j&\sum_{k=0}^i (\bij_{k,l,n+1}\nu_{n+1}(k,l)-\bijkln\nu_{n}(k,l))\\
= &\sum_{l=0}^j\sum_{k=0}^i \left[\bijklN\left((1-\frac{\delta_{kl}}{n})\nukln+c_1(k-1+\lambda)\frac{\nu_n(k-1,l)}{n}\right.\right.\\
&   \left.\left.+c_2(l-1+\mu)\frac{\nu_n(k,l-1)}{n}+\alpha\ind_{\{(i,j)=(0,1)\}}+\gamma \ind_{\{(i,j)=(1,0)\}}\right)-\bijkln\nukln\right]\\
\intertext{\sid{and identifying summands corresponding to
  $(k,l)=(i,j),\,(k,l)=(i-1,j),\,(k,l)=(i,j-1)$ and then the rest down
  to $(k,l)=(0,1), \, (k,l)=(1,0)$ yields}}
=& \bij_{i,j,n+1}\left(\nuij(1-\frac{\delta_{ij}}{n})+c_1(i-1+\lambda)\frac{\nu_n(i-1,j)}{n}+c_2(j-1+\mu)\frac{\nu_n(i,j-1)}{n}\right)\\
&\qquad -\bij_{i,j,n}\nuij\\
&+ \bij_{i-1,j,n+1}\left(\nu_n(i-1,j)(1-\frac{\delta_{i-1,j}}{n})+c_1(i-2+\lambda)\frac{\nu_n(i-2,j)}{n}+c_2(j-1+\mu)\frac{\nu_n(i-1,j-1)}{n}\right)\\
&\qquad -\bij_{i-1,j,n}\nu_n(i-1,j)\\
&+ \bij_{i,j-1,n+1}\left(\nu_n(i,j-1)(1-\frac{\delta_{i,j-1}}{n})+c_1(i-1+\lambda)\frac{\nu_n(i-1,j-1)}{n}+c_2(j-2+\mu)\frac{\nu_n(i,j-2)}{n}\right)\\
&\qquad -\bij_{i,j-1,n}\nu_n(i,j-1)\\
&+\ldots 
+ \bij_{0,1,n+1}\left(\alpha+\nu_n(0,1)(1-\frac{\delta_{01}}{n})\right)-\bij_{0,1,n}\nu_n(0,1)\\
&+ \bij_{1,0,n+1}\left(\gamma+\nu_n(1,0)(1-\frac{\delta_{10}}{n})\right)-\bij_{1,0,n}\nu_n(1,0)\\
=&\alpha\bij_{0,1,n+1}+\gamma\bij_{1,0,n+1}+\Biggl[\underbrace{
   \nuij\left(\bij_{i,j,n+1}\bigl(1-\frac{\delta_{ij}}{n}\bigr)-\bij_{i,j,n}\right)}_{=0
   \text{ by }\eqref{rec1}}\\
&+\sum_{k=0}^{i-1}\nu_n(k,j)\underbrace{
\left(\bij_{k,j,n+1}\bigl(1-\frac{\delta_{kj}}{n}\bigr)+\bij_{k+1,j,n+1}\frac{c_1(k+\lambda)}{n}-\bij_{k,j,n}\right)
}_{=0
  \text{ by }\eqref{rec2}}\\
&+ \sum_{l=0}^{j-1}\nu_n(i,l)\underbrace{
\left(\bij_{i,l,n+1}\bigl(1-\frac{\delta_{il}}{n}\bigr)+\bij_{i,l+1,n+1}\frac{c_2(l+\mu)}{n}-\bij_{i,l,n}\right)}_{=0
  \text{ by }\eqref{rec3}
  }\\
&+
  \sum_{l=0}^{j-1}\sum_{k=0}^{i-1}\underbrace{
\left(\bijklN\bigl(1-\frac{\delta_{kl}}{n}\bigr)+\bij_{k+1,l,n+1}\frac{c_1(k+\lambda)}{n}+\bij_{k,l+1,n+1}\frac{c_2(l+\mu)}{n}-\bij_{k,l,n}\right)}_{=0
  \text{ by }\eqref{rec4}}
  \Biggr]\\
=&\alpha\bij_{0,1,n+1}+\gamma\bij_{1,0,n+1}.
\end{align*}
So \eqref{mgdiff} now becomes
\be\label{mgdiffsimp}
M_{n+1}(i,j)-M_{n}(i,j)=
\sum_{l=0}^j\sum_{k=0}^i \Bigl(\bij_{k,l,n+1}N_{n+1}(k,l)-\bijkln N_{n}(k,l)\Bigr)
-(\alpha\bij_{0,1,n+1}+\gamma\bij_{1,0,n+1}).
\ee

\subsubsection{Conditional covariances.}\label{conditcov}
In order to use the multivariate martingale central limit theorem as specified in Proposition~\ref{mgclt}, we need to calculate the asymptotic form of the following quantity:
\be\label{scalediff}
\EE\left[\left(\frac{M_{n+1}(i,j)-M_{n}(i,j)}{\prod_{d=1}^n \left(1-\frac{\delta_{ij}}{d}\right)^{-1}}\right)\left(\frac{M_{n+1}(s,t)-M_{n}(s,t)}{\prod_{d=1}^n \left(1-\frac{\delta_{st}}{d}\right)^{-1}}\right)\middle|\FFF_n\right],
\ee for fixed pairs $(i,j)$ and $(s,t)$.
From \eqref{mgdiffsimp} we know that we need to consider in particular
\begin{align}
\bij_{k,l,n+1} & N_{n+1}(k,l)-\bijkln N_{n}(k,l)\nonumber \\
=& \bij_{k,l,n+1}N_{n+1}(k,l)\nonumber \\
 &-
   \left(\bijklN\left(1-\frac{\delta_{kl}}{n}\right)+\bij_{k+1,l,n+1}\frac{c_1(k+\lambda)}{n}+\bij_{k,l+1,n+1}\frac{c_2(l+\mu)}{n}\right)N_n(k,l)\nonumber \\  
\intertext{\sid{(where we applied \eqref{rec4})}}
=& \bij_{k,l,n+1}(N_{n+1}(k,l)-N_n(k,l))\nonumber \\  
 &+\frac{N_n(k,l)}{n}\left(\delta_{kl}+\bijklN-c_1(k+\lambda)\bij_{k+1,l,n+1}-c_2(l+\mu)\bij_{k,l+1,n+1}\right).\label{extra}
\end{align}
Recall \eqref{convij}  gives $N_n(k,l)/n\rightarrow p_{kl}$
a.s. as $n\rightarrow\infty$. So dealing with \eqref{scalediff} means 
we must calculate the 
asymptotic form of \sid{the conditional moments of} 
$$\sid{\Delta_{n+1}(i,j):=N_{n+1}(i,j)-\Nij.}$$

Observe that 
\begin{align}
\Delta_{n+1}(0,1)&=
\begin{cases}
1 & \text{w.p. } \alpha,\\
-1 & \text{w.p. } \delta_{01}\frac{N_n(0,1)}{n},\\
0 & \text{otherwise};
\end{cases}\label{D01}\\
\Delta_{n+1}(1,0)&=
\begin{cases}
1 & \text{w.p. } \gamma,\\
-1 & \text{w.p. } \delta_{10}\frac{N_n(1,0)}{n},\\
0 & \text{otherwise}; 
\end{cases}\label{D10} \\
\Delta_{n+1}(k,l) &=
\begin{cases}
1 & \text{w.p. } c_1(k-1+\lambda)\frac{N_n(k-1,l)}{n}+c_2(l-1+\mu)\frac{N_n(k,l-1)}{n},\\
-1 & \text{w.p. } \delta_{kl}\frac{N_n(k,l)}{n},\\
0 & \text{otherwise},
\end{cases}\label{Dkl}\end{align}
for $(k,l)\notin\{(0,1),(1,0)\}$. \sid{For instance, to justify 
\eqref{D01}, we create a $(0,1)$-node when node
 $n+1$ is born and
attaches to $V_n$ but we destroy a $(0,1)$-node if either $n+1$ is
born and attaches to a $(0,1)$-node or $v\in V_n$ attaches to $n+1$
and has degree $(0,1)$.}
 Then using \eqref{convij}, \eqref{pij} and \eqref{D01}--\eqref{Dkl}, for each pair $(k,l)$, 
\be\label{Dklasy}\EE(\Delta_{n+1}(k,l)|\FFF_n)\rightarrow p_{kl},\quad\text{a.s. as }n\rightarrow\infty.\ee
Therefore, \sid{{from }\eqref{prod}}
\begin{align*}
\frac{M_{n+1}(i,j)-M_{n}(i,j)}{\prod_{d=1}^n \left(1-\frac{\delta_{ij}}{d}\right)^{-1}}
=&\frac{M_{n+1}(i,j)-M_{n}(i,j)}{b_{i,j,n+1}^{(i,j)}}
\\
\intertext{and applying \eqref{mgdiffsimp} and then \eqref{extra}, we
  have as $n\to\infty$,}
=&\sid{\sum_{l=0}^j\sum_{k=0}^i \Biggr[ } \frac{N_n(k,l)}{n}\left(\delta_{kl}\frac{\bijklN}{b_{i,j,n+1}^{(i,j)}}-c_1(k+\lambda)\frac{\bij_{k+1,l,n+1}}{b_{i,j,n+1}^{(i,j)}}-c_2(l+\mu)
\frac{\bij_{k,l+1,n+1}}{b_{i,j,n+1}^{(i,j)}}\right)\\
&\quad
  +\frac{\bij_{k,l,n+1}}{b_{i,j,n+1}^{(i,j)}}\Delta_{n+1}(k,l)\sid{\Biggr]}
-\sid{\frac{(\alpha\bij_{0,1,n+1}+\gamma\bij_{1,0,n+1})}{b_{i,j,n+1}^{(i,j)}}      } \\
{\sim}& \sid{\sum_{l=0}^j\sum_{k=0}^i }
        p_{kl}\left(\delta_{kl}\xiij_{kl}-c_1(k+\lambda)\xiij_{k+1,l}-c_2(l+\mu)\xiij_{k,l+1}\right)+\xiij_{kl}\Delta_{n+1}^{(k,l)}\\
&\quad -\sid{(\alpha   
        \xi_{01}^{(i,j)} +\gamma \xi_{10}^{(i,j)}), }
\end{align*}
according to \eqref{convij} and definition of $\xiij_{kl}$ given in
\eqref{xi}, \eqref{limbil} and \eqref{limbkl}.
%In order to use the multivariate martingale central limit theorem, we need to calculate the asymptotic form of the following quantity:
%\be\label{scalediff}
%\EE\left[\left(\frac{M_{n+1}(i,j)-M_{n}(i,j)}{\prod_{d=1}^n \left(1-\frac{\delta_{ij}}{d}\right)^{-1}}\right)\left(\frac{M_{n+1}(s,t)-M_{n}(s,t)}{\prod_{e=1}^n \left(1-\frac{\delta_{st}}{e}\right)^{-1}}\right)\middle|\FFF_n\right],
%\ee for fixed pairs $(i,j)$ and $(s,t)$. 

Recall \eqref{Dklasy}, we see that as $n\rightarrow\infty$ the
\sid{conditional expectation} in \eqref{scalediff} is equivalent to 
\begin{align}
\EE & \Biggl\{ \Biggl[\sum_{l=0}^j \sum_{k=0}^i
      \Bigl\{p_{kl}\Bigl(\delta_{kl}\xiij_{kl}-c_1(k+\lambda)\xiij_{k+1,l}-c_2(l+\mu)\xiij_{k,l+1}\Bigr)+\xiij_{kl}\Delta_{n+1}^{(k,l)}\Bigr\}-(\alpha\xiij_{01}+\gamma\xiij_{10})\Biggr]\nonumber\\   
& \times \Biggl[\sum_{f=0}^t\sum_{h=0}^s
        \left\{p_{hf}\left(\delta_{hf}\xist_{hf}-c_1(h+\lambda)\xist_{h+1,f}-c_2(f+\mu)\xiij_{h,f+1}\right)+\xist_{hf}\Delta_{n+1}^{(h,f)}\right\}\nonumber
  \\
&\qquad
  -(\alpha\xist_{01}+\gamma\xist_{10})\Biggr] \Bigg|
  \FFF_n\Biggr\}\nonumber\\    
\intertext{\sid{and evaluating the product as four terms gives}}
&\sim \left[\sum_{l=0}^j\sum_{k=0}^i p_{kl}(\delta_{kl}\xiij_{kl}-c_1(k+\lambda)\xiij_{k+1,l}-c_2(l+\mu)\xiij_{k,l+1})-(\alpha\xiij_{01}+\gamma\xiij_{10})\right]\nonumber\\
     &\times \left[\sum_{f=0}^t\sum_{h=0}^s p_{hf}(\delta_{hf}\xist_{hf}-c_1(h+\lambda)\xist_{h+1,f}-c_2(f+\mu)\xiij_{h,f+1}-(\alpha\xist_{01}+\gamma\xist_{10})\right]\nonumber\\
&+ \left[\sum_{l=0}^j\sum_{k=0}^i p_{kl}(\delta_{kl}\xiij_{kl}-c_1(k+\lambda)\xiij_{k+1,l}-c_2(l+\mu)\xiij_{k,l+1})-(\alpha\xiij_{01}+\gamma\xiij_{10})\right]\nonumber\\
&\qquad  \times\left(\sum_{f=0}^t\sum_{h=0}^s \xist_{hf}p_{hf}\right)\nonumber\\
&+  \left[\sum_{f=0}^t\sum_{h=0}^s
  p_{hf}(\delta_{hf}\xist_{hf}-c_1(h+\lambda)\xist_{h+1,f}-c_2(f+\mu)\xiij_{h,f+1})-(\alpha\xist_{01}+\gamma\xist_{10})\right]\nonumber \\ 
&\qquad  \times\left(\sum_{l=0}^j\sum_{k=0}^i \xiij_{kl}p_{kl}\right)\nonumber\\
&+ \EE\left[\left(\sum_{l=0}^j\sum_{k=0}^i \xiij_{kl}\DklN\right)
                \left(\sum_{f=0}^t\sum_{h=0}^s \xist_{hf}\DhfN\right)\middle|\FFF_n\right]\nonumber\\
&=: A(i,j,s,t)+\EE\left[\left(\sum_{l=0}^j\sum_{k=0}^i \xiij_{kl}\DklN\right)
                \left(\sum_{f=0}^t\sum_{h=0}^s
  \xist_{hf}\DhfN\right)\middle|\FFF_n\right]                \label{eqsdiff}\\
&=A(i,j,s,t)+\sum_{(k,l)}\sum_{(h,f)}\xiij_{kl}\xist_{hf}\EE\left[\DklN\DhfN
  | \FFF_n \right]. \nonumber 
\end{align}

Therefore, we need the asymptotic form of the sum
\begin{align}
\sum_{(k,l)}\sum_{(h,f)}\xiij_{kl}\xist_{hf}\EE\left[\DklN\DhfN\middle|\FFF_n\right], \label{sumcross}                
\end{align}
and we  divide the summation in \eqref{sumcross} into four different cases.

\underline{Case I}: 
With probability $c_1(r-1+\lambda)\frac{N_n(r-1,q)}{n}$, 
a new edge from $n+1$ to some existing node $v\in V_n$ with
$D_n(v)=(r-1,q)$ is created \sid{and this necessitates
\begin{align*}
&\Delta_{n+1}(r-1,q)=-1,&&\quad \text{since an $(r-1,q)$-node is destroyed,}\\
&\Delta_{n+1}(r,q)=1,&&\quad \text{since an $(r,q)$-node is created,}\\
&\Delta_{n+1}(0,1)=1,&&\quad \text{since a $(0,1)$-node is created.}
\end{align*}
%Denote this case by event $E_1$, then $\PP(E_1)=c_1(r-1+\lambda)\frac{N_n(r-1,q)}{n}$.
%On $E_1$ we have
%\begin{align}
%\sum_{(k,l)}\sum_{(h,f)} & \xiij_{kl}\xist_{hf}\EE\left[\DklN\DhfN\middle|\FFF_n;E_1\right]\nonumber\\
%&=\sum_{q=0}^n\sum_{r=0}^n\left(\xiij_{rq}+\xiij_{01}-\xiij_{r-1,q}\right)\left(\xist_{rq}+\xist_{01}-\xist_{r-1,q}\right). \label{E1}
%\end{align}
The other cases follow similar reasoning:}

\underline{Case II}: 
With probability $c_2(q-1+\mu)\frac{N_n(r,q-1)}{n}$, 
a new edge from some existing node $v\in V_n$ (with $D_n(v)=(r,q-1)$) to $n+1$ is created such that
\[\Delta_{n+1}(r,q-1)=-1,\quad \Delta_{n+1}(r,q)=1,\quad \Delta_{n+1}(1,0)=1.\]
%On $E_2$,
%\begin{align}
%\sum_{(k,l)}\sum_{(h,f)} & \xiij_{kl}\xist_{hf}\EE\left[\DklN\DhfN\middle|\FFF_n;E_2\right]\nonumber\\
%&= \left(\xiij_{rq}+\xiij_{01}-\xiij_{r,q-1}\right)\left(\xist_{rq}+\xist_{01}-\xist_{r,q-1}\right),\label{E2}
%\end{align}
%with $P(E_2)=c_2(q-1+\mu)\frac{N_n(r,q-1)}{n}$.\\

\underline{Case III}: With probability $c_1(r+\lambda)\frac{N_n(r,q)}{n}$, a new edge from $n+1$ to some existing node $v\in V_n$ (with $D_n(v)=(r,q)$) is created such that
\[\Delta_{n+1}(r,q)=-1,\quad \Delta_{n+1}(r+1,q)=1,\quad \Delta_{n+1}(0,1)=1.\]

\underline{Case IV}: With probability $c_2(q+\mu)\frac{N_n(r,q)}{n}$, a new edge from some existing node $v\in V_n$ (with $D_n(v)=(r,q)$) to $n+1$ is created such that
\[\Delta_{n+1}(r,q)=-1,\quad \Delta_{n+1}(r,q+1)=1,\quad \Delta_{n+1}(1,0)=1.\]

Take Case I as an example, we see that
\begin{equation}\label{prodelta}
\DklN\DhfN=
 \begin{cases}
 1  & \text{ if }  ((k,l),(h,f))\in\{((r,q),(r,q)),((r-1,q), (r-1,q)), \\
 &   ((0,1),(0,1)),((r,q),(0,1)),((0,1),(r,q))\};\\
 -1 & \text{ if }\left((k,l),(h,f)\right)\in\{((r-1,q),(r,q)),((r-1,q),(0,1)),\\
 & ((r,q),(r-1,q)),((0,1),(r-1,q))\};\\
 0 & \text{ otherwise.}
 \end{cases}
\end{equation}
Let $E_1^{(r,q)}$ denote the event described in Case I \sid{where node
  $n+1$ attaches to $v\in V_n$ with $D_n(v)=(r-1,q)$.}
Then on the event $\mathcal{E}_1:=\bigcup_{(r,q)}E_1^{(r,q)}$,
\eqref{prodelta} gives \sid{asymptotically}
\begin{align}
\sum_{(k,l)}\sum_{(h,f)} &
            \xiij_{kl}\xist_{hf}\EE\left[\DklN\DhfN \sid{1_{\mathcal{E}_1}}\middle|\FFF_n\right]\nonumber\\     
&= \sum_{(r,q)}\PP(E_1^{(r,q)})\left(\xiij_{rq}+\xiij_{01}-\xiij_{r-1,q}\right)\left(\xist_{rq}+\xist_{01}-\xist_{r-1,q}\right)\nonumber\\
&=\sum_{q=0}^n\sum_{r=0}^n c_1(r-1+\lambda)\frac{N_n(r-1,q)}{n}\left(\xiij_{rq}+\xiij_{01}-\xiij_{r-1,q}\right)\left(\xist_{rq}+\xist_{01}-\xist_{r-1,q}\right). \label{case1}
\end{align}

Define $\mathcal{E}_2$, $\mathcal{E}_3$ and $\mathcal{E}_4$ in the
same way with respect to Case II, III and IV, \sid{and} then similar
calculations to \eqref{case1} give $\sum_{(k,l)}\sum_{(h,f)}
\xiij_{kl}\xist_{hf}\EE\left[\DklN\DhfN\middle|\FFF_n;\mathcal{E}_i\right]$
for $i=2,3,4$. Also, \eqref{D01} and \eqref{D10} show that
$\EE\left[\left(\Delta_{n+1}(0,1)\right)^2\middle|\FFF_n\right]$ and
$\EE\left[\left(\Delta_{n+1}(1,0)\right)^2\middle|\FFF_n\right]$ take
different forms from the other cases (cf. \eqref{Dkl}), so we still
need to compensate for this. 

Considering the case where $(k,l)=(h,f)=(0,1)$, we have, by \eqref{D01},
\begin{align*}
\xiij_{01}\xist_{01}\EE\left[\left(\Delta_{n+1}(0,1)\right)^2\middle|\FFF_n\right]
&= \xiij_{01}\xist_{01}\left(\alpha+\delta_{01}\frac{N_n(0,1)}{n}\right)\\
&= \xiij_{01}\xist_{01}\left(\alpha+c_1\lambda\frac{N_n(0,1)}{n}+c_2(1+\mu)\frac{N_n(0,1)}{n}\right).
\end{align*}
Note that $\xiij_{01}\xist_{01}c_2(1+\mu)\frac{N_n(0,1)}{n}$ has been covered while calculating \eqref{sumcross} with respect to $\mathcal{E}_4$, so we only need to add $\xiij_{01}\xist_{01}\left(\alpha+c_1\lambda\frac{N_n(0,1)}{n}\right)$ to our computation. Similar arguments also apply to $(k,l)=(h,f)=(1,0)$, but instead we add $\left(\gamma+c_2\mu\frac{N_n(1,0)}{n}\right)\xiij_{10}\xist_{10}$ for compensation. 

Taking all these into account, we get
\begin{eqnarray*}
\lefteqn{\sum_{(k,l)}\sum_{(h,f)}\xiij_{kl}\xist_{hf}\EE\left[\DklN\DhfN\middle|\FFF_n\right]}\\
&=& \left(\alpha+c_1\lambda\frac{N_n(0,1)}{n}\right)\xiij_{01}\xist_{01}+\left(\gamma+c_2\mu\frac{N_n(1,0)}{n}\right)\xiij_{10}\xist_{10}\\
&&\:+\sum_{q=0}^n\sum_{r=0}^n\left\{
c_1(r-1+\lambda)\frac{N_n(r-1,q)}{n}\left(\xiij_{rq}+\xiij_{01}-\xiij_{r-1,q}\right)\left(\xist_{rq}+\xist_{01}-\xist_{r-1,q}\right)\right.\\
&&\:\left.+c_2(q-1+\mu)\frac{N_n(r,q-1)}{n}\left(\xiij_{rq}+\xiij_{01}-\xiij_{r,q-1}\right)\left(\xist_{rq}+\xist_{01}-\xist_{r,q-1}\right)\right.\\
&&\:\left.+c_1(r+\lambda)\frac{N_n(r,q)}{n}\left(\xiij_{r+1,q}+\xiij_{01}-\xiij_{rq}\right)\left(\xist_{r+1,q}+\xist_{01}-\xist_{rq}\right)\right.\\
&&\: \left.+c_2(q+\mu)\frac{N_n(r,q)}{n}\left(\xiij_{r,q+1}+\xiij_{01}-\xiij_{rq}\right)\left(\xist_{r,q+1}+\xist_{01}-\xist_{rq}\right)
\right\}.
\end{eqnarray*}
Here we also adopt the convention that $\xi_{kl}^{(i,j)}=0$ if either
$k>i$ or $l>j$, and that $N_n(r,q)=0$ whenever either both $r$ and $q$
are 0 or one of them is $-1$.

%Meanwhile, \eqref{D01}--\eqref{Dklasy} give that as $n\rightarrow\infty$,
%\begin{align*}
%\EE\left[\left(\Delta_{n+1}(0,1)\right)^2\middle|\FFF_n\right] &\rightarrow \alpha+\delta_{01}p_{01},\\
%\EE\left[\left(\Delta_{n+1}(1,0)\right)^2\middle|\FFF_n\right] &\rightarrow \gamma+\delta_{10}p_{10},\\
%\EE\left[\left(\DklN\right)^2\middle|\FFF_n\right] &\rightarrow c_1(k-1+\lambda)p_{k-1,l}+c_2(l-1+\mu)p_{k,l-1}+\delta_{kl}p_{kl},\\
%   \text{for }(k,l)&\notin \{(0,1),(1,0)\}.
%\end{align*}
%Thus,

%Hence, using the convention that $\xi_{kl}^{(i,j)}=0$ if either $k>i$ or $l>j$, 
Now applying \eqref{convij} again, we write
\begin{eqnarray}\label{croasy}
\lefteqn{\sum_{(k,l)}\sum_{(h,f)}\xiij_{kl}\xist_{hf}\EE\left[\DklN\DhfN\middle|\FFF_n\right]}\nonumber\\
&\rightarrow & \left(\alpha+c_1\lambda p_{01}\right)\xiij_{01}\xist_{01}+\left(\gamma+c_2\mu p_{10}\right)\xiij_{10}\xist_{10}\nonumber\\
&&\:+\sum_{q=0}^\infty\sum_{r=0}^\infty\left\{
c_1(r-1+\lambda) p_{r-1,q}\left(\xiij_{rq}+\xiij_{01}-\xiij_{r-1,q}\right)\left(\xist_{rq}+\xist_{01}-\xist_{r-1,q}\right)\right.\nonumber\\
&&\:\left.+c_2(q-1+\mu) p_{r,q-1}\left(\xiij_{rq}+\xiij_{01}-\xiij_{r,q-1}\right)\left(\xist_{rq}+\xist_{01}-\xist_{r,q-1}\right)\right.\nonumber\\
&&\:\left.+c_1(r+\lambda) p_{rq}\left(\xiij_{r+1,q}+\xiij_{01}-\xiij_{rq}\right)\left(\xist_{r+1,q}+\xist_{01}-\xist_{rq}\right)\right.\nonumber\\
&&\: \left.+c_2(q+\mu) p_{rq}\left(\xiij_{r,q+1}+\xiij_{01}-\xiij_{rq}\right)\left(\xist_{r,q+1}+\xist_{01}-\xist_{rq}\right)
\right\}=: B(i,j,s,t)
\end{eqnarray}
a.s. as $n\rightarrow\infty$. Putting \eqref{eqsdiff} and \eqref{croasy} together, we conclude that, with probability 1, 
\be\label{limijst}\EE\left[\left(\frac{M_{n+1}(i,j)-M_{n}(i,j)}{\prod_{d=1}^n \left(1-\frac{\delta_{ij}}{d}\right)^{-1}}\right)\left(\frac{M_{n+1}(s,t)-M_{n}(s,t)}{\prod_{d=1}^n \left(1-\frac{\delta_{st}}{d}\right)^{-1}}\right)\middle|\FFF_n\right]\rightarrow C(i,j,s,t),\ee
where $C(i,j,s,t):=A(i,j,s,t)+B(i,j,s,t)$.

Recall that $\bij_{i,j,n+1}=\prod_{d=1}^n(1-\delta_{ij}/d)^{-1}$. \sid{By} Stirling's formula, as $n\rightarrow\infty$
\be\label{aij}
\bij_{i,j,n+1}=\prod_{d=1}^n\frac{d}{d-\delta_{ij}}=\frac{\Gamma(n+1)}{\Gamma(n+1-\delta_{ij})\Gamma(1-\delta_{ij})}\sim \frac{n^{\delta_{ij}}}{\Gamma(1-\delta_{ij})},
\ee
so that as a function of $n$, $\bij_{i,j,n+1}$ is regularly varying with index $\delta_{ij}$. Therefore, \eqref{limijst} becomes
\begin{align}\label{tauij}
\EE\left[\left(\frac{M_{n+1}(i,j)-M_n(i,j)}{n^{\delta_{ij}}}\right)\left(\frac{M_{n+1}(s,t)-M_n(s,t)}{n^{\delta_{st}}}\right)\middle|\FFF_n\right] \rightarrow & \frac{C(i,j,s,t)}{\Gamma(1-\delta_{ij})\Gamma(1-\delta_{st})}\nonumber\\
&=: \tau(i,j,s,t).
\end{align}

\subsubsection{Applying the martingale central
  limit theorem.}
\sid{We now have the material necessary to verify the conditions in
Proposition \ref{mgclt}.}
\sid{Fix non-negative integers $I,O \in \{0,1,2,\ldots\}$ and  define
for $I\vee O+1\leq m\leq n$,}
$$
X_{n,m,i,j}= \frac{M_m(i,j)-M_{m-1}(i,j)}{n^{\delta_{ij}+1/2}}
               \;0\leq i\leq I, 0\leq j \leq O,
$$
\sid{and with   $(s,t)$ satisfying   $  0\leq s\leq I, \,0\leq t \leq
  O,$ also define}
\begin{align*}
G_{n,m}(i,j,s,t)&:=\EE(X_{n,m,i,j}X_{n,m,s,t}|\FFF_{m-1})\\
&= n^{-(\delta_{ij}+\delta_{st}+1)}\EE[(M_m(i,j)-M_{m-1}(i,j))(M_m(s,t)-M_{m-1}(s,t))|\FFF_{m-1}].
\end{align*}
We know from \eqref{tauij} that
\be\label{Gnasy}
nG_{n,n}(i,j,s,t)\rightarrow\tau(i,j,s,t),\quad \mbox{as }n\rightarrow\infty,
\ee
and that
\[G_{n,m}(i,j,s,t)=\frac{m^{\delta_{ij}+\delta_{st}}}{n^{1+\delta_{ij}+\delta_{st}}}mG_{m,m}(i,j,s,t).\]
Hence, by Karamata's theorem on integration of regularly varying functions, using \eqref{Gnasy} we have
\begin{align}\label{asyvar}
V_n(i,j,s,t)\sid{:}&= \sum_{m=I\vee O+1}^n G_{n,m}(i,j,s,t)=\frac{\sum_{m=I\vee O}^n m^{1+\delta_{ij}+\delta_{st}}G_{m,m}(i,j,s,t)}{n^{1+\delta_{ij}+\delta_{st}}}\nonumber\\
&\sim  \frac{n\cdot
  n^{1+\delta_{ij}+\delta_{st}}G_{n,n}(i,j,s,t)}{(1+\delta_{ij}+\delta_{st})n^{1+\delta_{ij}+\delta_{st}}}\sim
  \frac{\tau(i,j,s,t)}{1+\delta_{ij}+\delta_{st}} 
=\sigma^2(i,j,s,t).\end{align}
\sid{So the $\bigl(  (I+1)\times (O+1)\times (I+1)\times (O+1)\bigr)$
dimensional matrix \sidney{converges}
\begin{equation}\label{matrixAsyVar}
\Bigl(V_n(i,j,s,t); 0\leq i,s \leq I,\,0\leq j,t \leq O\Bigr) \to
\Sigma_{IO} =(\sigma^2(i,j,s,t))\end{equation}
as required by Propositon \ref{mgclt}.}
For each pair $(i,j)$ such that $0\leq i\leq I,\,0\leq j \leq O$, from
the definition of $M_n(i,j)$ in \eqref{mg}, 
\begin{align*}
\frac{M_n(i,j)}{n^{\delta_{ij}+1/2}}
&= \sum_{l=0}^j\sum_{k=0}^i
  \frac{b_{k,l,n}^{(i,j)}}{n^{\delta_{ij}}}\left(\frac{N_n(k,l)-\nu_n(k,l)}{\sqrt{n}}\right)\\ 
& \sim  \Gamma(1-\delta_{ij})\sum_{l=0}^j\sum_{k=0}^i\xi_{k,l}^{(i,j)}\left(\frac{N_n(k,l)-\nu_n(k,l)}{\sqrt{n}}\right)
\end{align*}
\sid{and this lets us write the matrix equation (with $o_p(1)$ terms dropped)
\begin{align}\label{KIO}
\Bigl( \frac{M_n(i,j)}{n^{\delta_{ij}+1/2}}; 0\leq i\leq I,0\leq j\leq O\Bigr)
&=:  K_{IO}\left(\frac{N_n(k,l)-\nu_n(k,l)}{\sqrt{n}}:\,  0\le k\le I,\, 0\le l\le O\right)^T,
\end{align}
}
\sid{where we think of $\bigl((N_n(k,l)-\nu_n(k,l)/{\sqrt{n}}, 0\le
k\le I,\, 0\le l\le O\bigr)$ as a $(I+1)\times(O+1)$ dimensional column
vector.}
\sid{Relation} \eqref{KIO} results from definitions of $\xi_{k,l}^{(i,j)}=\lim_{n\rightarrow\infty}b_{k,l,n+1}^{(i,j)}/b_{i,j,n+1}^{(i,j)}$:
\begin{align*}
 \frac{b_{k,l,n}^{(i,j)}}{b_{i,j,n}^{(i,j)}}
 &= \frac{b_{k,l,n}^{(i,j)}}{b_{i,j,n+1}^{(i,j)}}\frac{b_{i,j,n+1}^{(i,j)}}{b_{i,j,n}^{(i,j)}}\\
 &\stackrel{\eqref{rec1}}{=} \frac{b_{k,l,n}^{(i,j)}}{b_{i,j,n+1}^{(i,j)}} \left(1-\frac{\delta_{ij}}{n}\right)^{-1}\\
 &\stackrel{\eqref{rec4}}{=} \left\{\frac{b_{k,l,n+1}^{(i,j)}}{b_{i,j,n+1}^{(i,j)}}\left(1-\frac{\delta_{kl}}{n}\right)
 +\frac{b_{k+1,l,n+1}^{(i,j)}}{b_{i,j,n+1}^{(i,j)}}\frac{c_1(k+\lambda)}{n}
 +\frac{b_{k,l+1,n+1}^{(i,j)}}{b_{i,j,n+1}^{(i,j)}}\frac{c_2(l+\mu)}{n}\right\}\left(1-\frac{\delta_{ij}}{n}\right)^{-1}\\
 &\rightarrow \xi_{k,l}^{(i,j)},\quad\text{as n}\rightarrow \infty,
\end{align*}
provided that we set $b_{k,l,n}^{(i,j)}=0$ if either $k=i+1$ or $l=j+1$. Then similar to \eqref{aij},
\[b^{(i,j)}_{i,j,n}=\prod_{d=1}^{n-1}\frac{d}{d-\delta_{ij}}=\frac{\Gamma(n)}{\Gamma(n-\delta_{ij})\Gamma(1-\delta_{ij})}\sim \frac{n^{\delta_{ij}}}{\Gamma(1-\delta_{ij})},
\] thus giving the equivalence relationship in \eqref{KIO}.

\sidney{In order to  apply} Proposition~\ref{mgclt} to conclude
\eqref{clt}, we must verify 
\sidney{conditions (i) and (ii) of the Proposition.}
Condition  (i) of Proposition \ref{mgclt} is already satisfied by \eqref{asyvar}, so we just need to consider condition (ii).
Since by \eqref{D01}--\eqref{Dkl} the differences are bounded, i.e.
\[|(\Nij-\nuij)-(N_{n-1}(i,j)-\nu_{n-1}(i,j))|\le 2\mbox{ for all }(i,j),\]
then we claim that for $n$ large enough, the events $\{|X_{n,m,i,j}|>\varepsilon\}$ vanish for all $m\le n$ and all $(i,j)$. This can be observed from the following. For some constant $\kappa_{ij}$,
\bean
\{|X_{n,m,i,j}|>\varepsilon\}&=& \{|M_m(i,j)-M_{m-1}(i,j)|>\varepsilon n^{\delta_{i,j}+1/2}\}\\
&\subseteq & \{|\kappa_{ij}|m^{\delta_{ij}}>\varepsilon n^{\delta_{ij}+1/2}\}\\
&\subseteq & \{|\kappa_{ij}|n^{\delta_{ij}}>\varepsilon n^{\delta_{ij}+1/2}\}
\eean
So with probability converging to 1 as $n\rightarrow\infty$, the indicator functions $\ind_{\{|X_{n,m,i,j}|>\varepsilon\}}$ vanish. This verifies the second condition, i.e.

Recall that calculations in \eqref{asyvar} \sid{and
\eqref{matrixAsyVar} } gives the covariance 
matrix $\Sigma_{IO}$.  \sidney{A}pplying Proposition~\ref{mgclt} yields
\be\label{concl1}
K_{IO}\left(\frac{N_n(i,j)-\nu_n(i,j)}{\sqrt{n}}:\,  0\le i\le I,\, 0\le j\le O\right)^T\Rightarrow N(0,\Sigma_{I,O})
\ee
in $\RR^{(I+1)(O+1)}$. \sid{If we a}ssume further that $K_{IO}$ is invertible, then the convergence in \eqref{concl1} can be rewritten as
\[
\left(\frac{N_n(i,j)-\nu_n(i,j)}{\sqrt{n}}:\,  0\le i\le I,\, 0\le j\le O\right)\Rightarrow N(0,K^{-1}_{IO}\Sigma_{I,O}K^{-T}_{IO}).
\]
Applying Lemma~\ref{degseq}, we can then obtain that for fixed $I$ and $O$, \eqref{thm} holds.
%\[\left(\sqrt{n}\left(\frac{N_n(i,j)}{n}-p_{ij}\right):\,  0\le i\le I,\, 0\le j\le O\right)\Rightarrow N(0,K^{-1}_{IO}\Sigma_{I,O}K^{-T}_{IO}).\]

%We are now left with checking the two conditions in
%Proposition~\ref{mgclt}. 
%\sid{??Are you doing this here? What does the
%next paragraph verify? Give it a name or number and refer to it. Isn't
%the first condition  (i) of Proposition \ref{mgclt} verified by \eqref{asyvar}?}

To avoid non-degenerate limits, we need to make sure that the asymptotic variances given 
in matrix $\Sigma_{I,O}$ are positive for fixed $I$ and $O$. It suffices to check that 
for $0\le i \le I$, $0\le j\le O$,
\be\label{posvar}
\lim_{n\rightarrow\infty}\frac{\VV_n(i,j)}{n}: =\lim_{n\to\infty} \frac{\Var(\Nij)}{n}> 0.
\ee
From the definition,
\begin{align*}
\VV_{n+1}(i,j)&= \EE\left[\left(\NijN\right)^2\right]-\left(\nuijN\right)^2\\
&=\EE\left[\EE\left((\Nij+\DijN)^2\middle|\FFF_n\right)\right]-\left(\nuijN\right)^2.
%&&\:-\left(\nuij(1-\frac{\delta_{ij}}{n})+c_1(i-1+\lambda)\frac{\nu_n(i-1,j)}{n}
%+c_2(j-1+\mu)\frac{\nu_n(i,j-1)}{n}\right)^2
\end{align*}
For $(i,j)\notin\{(0,1),(1,0)\}$, we have from \eqref{Dkl} 
\begin{align}\label{var1}
\EE & \left((\Nij+\DijN)^2\middle|\FFF_n\right)\nonumber\\
=& (\Nij)^2+2\Nij\EE(\Delta_{n+1}(i,j)|\FFF_n)
+\EE\left((\DijN)^2\middle|\FFF_n\right)\nonumber\\
=& (\Nij)^2+2\Nij\left(c_1(i-1+\lambda)\frac{N_n(i-1,j)}{n}+c_2(j-1+\mu)\frac{N_n(i,j-1)}{n}-\delta_{ij}\frac{\Nij}{n}\right)\nonumber\\
&+  c_1(i-1+\lambda)\frac{N_n(i-1,j)}{n}+c_2(j-1+\mu)\frac{N_n(i,j-1)}{n}+\delta_{ij}\frac{\Nij}{n},
\end{align}
and using~\eqref{nudef}
\begin{align}\label{var2}
 (& \nuijN)^2 \nonumber\\
=& \left(\nuij+c_1(i-1+\lambda)\frac{\nu_n(i-1,j)}{n}
+c_2(j-1+\mu)\frac{\nu_n(i,j-1)}{n}-\delta_{ij}\frac{\nuij}{n}\right)^2\nonumber\\
=& (\nuij)^2+2\nuij\left[c_1(i-1+\lambda)\frac{\nu_n(i-1,j)}{n}
+c_2(j-1+\mu)\frac{\nu_n(i,j-1)}{n}-\delta_{ij}\frac{\nuij}{n}\right]\nonumber\\
&+ \left[c_1(i-1+\lambda)\frac{\nu_n(i-1,j)}{n}
+c_2(j-1+\mu)\frac{\nu_n(i,j-1)}{n}-\delta_{ij}\frac{\nuij}{n}\right]^2.
\end{align}
Therefore, taking the expectation on both sides of \eqref{var1} and subtracting \eqref{var2} from it give
\bean
\VV_{n+1}(i,j)&=& \VV_n(i,j)\left(1-\frac{2\delta_{ij}}{n}\right)
+\frac{2c_1(i-1+\lambda)}{n}\EE[\Nij N_n(i-1,j)-\nuij\nu_n(i-1,j)]\\
&&\:+ \frac{2c_2(j-1+\mu)}{n}\EE[\Nij N_n(i,j-1)-\nuij\nu_n(i,j-1)]+R_{n+1}(i,j),
\eean
where as $n\rightarrow\infty$,
\bean
R_{n+1}(i,j)& \rightarrow & c_1(i-1+\lambda)p_{i-1,j}+c_2(j-1+\mu)p_{i,j-1}+\delta_{ij}p_{ij}\\
&&\:- [c_1(i-1+\lambda)p_{i-1,j}+c_2(j-1+\mu)p_{i,j-1}-\delta_{ij}p_{ij}]^2\\
&=& (1+2\delta_{ij})p_{ij}-p_{ij}^2\\
&=& 2\delta_{ij}p_{ij}+p_{ij}(1-p_{ij})> 0,
\eean
since $p_{ij}\in (0,1]$ and $\delta_{ij}> 0$. Note that here $p_{ij}\neq 0$ for all $(i,j)$: the recursion in \eqref{pij} shows that both $p_{01}, p_{10}>0$ as we assume $\alpha,\gamma>0$; it also follows that $p_{ij}=0$ for all $(i,j)\notin\{(0,1),(1,0)\}$ if we assume $p_{ij}=0$ for some $(i,j)$, which is impossible since we initiate the graph with a single node $v$ and $D_1(v)=(1,1)$.

Let $L_{ij}$ denote the limit of $R_{n+1}(i,j)$, then there exists $n_0$ such that for all $n\ge n_0$, $R_n(i,j)\ge\frac{1}{2}L_{ij}$. Also,
\begin{align*}
\EE[& \Nij N_n(i-1,j)-\nuij\nu_n(i-1,j)]\\
&=\mbox{cov}(\Nij, N_n(i-1,j))\ge -(V_n(i,j))^{1/2}(V_n(i-1,j))^{1/2},
\end{align*}
and similarly 
\[
\EE[\Nij N_n(i,j-1)-\nuij\nu_n(i,j-1)]\ge -(V_n(i,j))^{1/2}(V_n(i,j-1))^{1/2}.
\]
 We can now prove \eqref{posvar} by induction. The base case when $n=1$ is trivial. For $n\ge 2$, suppose that $\VV_n(i-1,j)\ge a_{i-1,j}n$ and $\VV_n(i,j-1)\ge a_{i,j-1}n$ for some $a_{i-1,j}, a_{i,j-1}>0$, then for all $n\ge n_0$,
\bean
\VV_{n+1}(i,j)&\ge & \VV_n(i,j)\left(1-\frac{\delta_{ij}}{n}\right)
-2c_1(i-1+\lambda)\left(\frac{\VV_n(i,j)}{n}\right)^{1/2}\left(\frac{\VV_n(i-1,j)}{n}\right)^{1/2}\\
&&\:-2c_2(j-1+\mu)\left(\frac{\VV_n(i,j)}{n}\right)^{1/2}\left(\frac{\VV_n(i,j-1)}{n}\right)^{1/2}
+\frac{1}{2}L_{ij}.
\eean
We can therefore conclude that 
\[\VV_{n+1}(i,j)\ge \VV_n(i,j)\left(1-\frac{K_1^{(i,j)}}{n}\right)
-K_2^{(i,j)}\left(\frac{\VV_n(i,j)}{n}\right)^{1/2}+\frac{1}{2}L_{ij} \quad\forall n\ge n_0,\]
where $K_1^{(i,j)},K_2^{(i,j)}>0$ positive constants.

If $\frac{\VV_n(i,j)}{n}\le \left(\frac{L_{ij}}{4K_2^{(i,j)}}\right)^2$, then
\[-K_2^{(i,j)}\left(\frac{\VV_n(i,j)}{n}\right)^{1/2}\ge\frac{1}{4}L_{ij}.\]
If $\frac{\VV_n(i,j)}{n}\ge \left(\frac{L_{ij}}{4K_2^{(i,j)}}\right)^2$, then
\[\left(\frac{\VV_n(i,j)}{n}\right)^{1/2}\le\frac{4K_2^{(i,j)}}{L_{ij}}\cdot \frac{\VV_n(i,j)}{n}.\]
In either case, for all $n\ge n_0$ we still have
\bean
\VV_{n+1}(i,j)&\ge & \VV_n(i,j)\left(1-\frac{K_1^{(i,j)}}{n}\right)
-K_2^{(i,j)}\frac{4K_2^{(i,j)}}{L_{ij}}\cdot \frac{\VV_n(i,j)}{n}+\frac{1}{4}L_{ij}\\
&=& \VV_n(i,j)(1-\frac{K_1^{(i,j)}+4\left(K_2^{(i,j)}\right)^2/L_{ij}}{n})+\frac{1}{4}L_{ij}.
\eean
Since $K^{(i,j)}:=K_1^{(i,j)}+4\left(K_2^{(i,j)}\right)^2/L_{ij}>0$, then for $n\ge n_0$,
\bean
\VV_{n+1}(i,j)&\ge & \frac{1}{4}L_{ij}\left[1+\left(1-\frac{K^{(i,j)}}{n}\right)+\left(1-\frac{K^{(i,j)}}{n}\right)^2+\cdots+\left(1-\frac{K^{(i,j)}}{n}\right)^{n-n_0}\right]\\
&=& \frac{1}{4}L_{ij}\frac{1-(1-K^{(i,j)}/n)^{n-n_0+1}}{K^{(i,j)}/n}\\
&\sim & \frac{L_{ij}}{4K^{(i,j)}}(1-e^{-K^{(i,j)}})n> 0, \mbox{    as }n\rightarrow\infty.
\eean
So we are done with proving \eqref{posvar}, thus completing the proof for Theorem~\ref{thm}.

\bibliographystyle{plain}
\bibliography{bib-directedCLT}

\bigskip\noindent
Tiandong Wang, School of Operations Research and Information Engineering, Cornell University, Ithaca NY 14853. 
Email: tw398@cornell.edu

\bigskip\noindent
Sidney I. Resnick, School of Operations Research and Information Engineering, Cornell University, Ithaca NY 14853. 
Email: sir1@cornell.edu

\end{document}